\magnification=1200

\def\C{{\cal C}}
\def\M{{\cal M}}
\def\R{{\cal R}}
\def\hal{{\vrule height 10pt width 4pt depth 0pt}}

\centerline{\bf Automatic Convexity}
\medskip

\centerline{Charles A. Akemann and Nik Weaver\footnote{*}{Second author
supported by NSF grant DMS-0070634}}
\bigskip
\bigskip

{\narrower{\it
In many cases the convexity of the image of a linear map with
range is $\R^n$ is automatic because of the facial structure of the
domain of the map. We develop a four step procedure for proving this kind
of ``automatic convexity''.  To make this procedure more efficient, we
prove two new theorems that identify the facial structure of the
intersection of a convex set with a subspace in terms of the
facial structure of the original set.

Let $K$ be a convex set in a real linear space $X$ and let $H$
be a subspace of X that meets $K$. In Part I we show that the faces of
$K\cap H$ have the form $F\cap H$ for a face $F$ of $K$.  Then we extend
our intersection theorem to the case where $X$ is a locally convex linear
topological space, $K$ and $H$ are closed, and $H$ has finite codimension
in $X$. In Part II we use our procedure to ``explain'' the convexity of
the numerical range (and some of its generalizations) of a complex matrix.
In Part III we use the topological version of our intersection theorem to
prove a version of Lyapunov's theorem with finitely many linear
constraints. We also extend Samet's continuous lifting theorem to the
same constrained siuation.
\bigskip}}

Historically there have been several theorems that concluded, unexpectedly,
even mysteriously at first, that a certain set in $\R^n$ is convex.  Perhaps
the two best known examples are the convexity of the numerical range of an
$n\times n$ complex matrix [Hau, T] and Lyapunov's theorem on the
convexity of the range of a vector measure [Ly].  In each of these cases
the set in question is the image under some apparently non-linear map of a
non-convex set. Each of these theorems has been generalized in many
directions.  Until  the work of Lindenstrauss [Li], Lyapunov's theorem
remained a mystery with several complicated, yet incomplete, proofs
(including Lyapunov's and a later proof by Halmos [Hal-1]) in the
literature. See [AA] for a discussion of Lyapunov's theorem and
generalizations. As for the convexity of the numerical range, while the
proofs in the literature have been complete, and they have gotten
steadily simpler, the mystery of the appearance of convexity has remained
(see [HJ, p. 78], [P] and [GR, sections 1.1 and 5.5]).

In [AA] a number of automatic convexity theorems related to Lyapunov's
Theorem were proved.  The key to those theorems is given in
[AA, Theorem 1.6 and Corollary 1.7], which we restate here, correcting
misprints, after introducing some notation.
\bigskip

\noindent {\bf NOTATION:}  $K$ denotes a convex set in a real linear space
$X$.  For any distinct points $x, y \in X$ let $(x,y)$ denote the line
segment joining $x$ and $y$, excluding the end points. $E(K)$ denotes the
set of extreme points of
$K$. If $K$ is not a singleton, the facial dimension [AA, p. 10] of
$K$ is defined to be $inf\{dim(Q) : Q $ is a nonsingleton face of $K\}.$
(Facial dimension $\infty$ is quite possible and especially interesting as
we shall see in Part III of this paper.)  For any subset $F$ of $K$ let
$G(K,F)$ denote the smallest face of  $K$ containing $F$.  In [AA] this
concept was defined and developed for singleton sets $F = \{v\}$, where the
notation $G(k,v)$ was used.
\bigskip

\noindent
{\bf [AA, 1.6].} If $K$
has facial dimension $>n$, $\Psi$ is an affine map of $K$ into $\R^n$,
and $v \in K$, then every extreme point of $\Psi ^{-1}(\Psi(v))$ is an
extreme point of $K$.
\medskip
\noindent
{\bf [AA, 1.7].}  If $X$ is a locally convex space, $K$ is compact with
facial dimension $> n$, and $\Psi$ is a continuous affine map from
$K$ into $\R^n$, then $\Psi(E(K)) = \Psi(K)$.
\medskip

The form of [AA, 1.7] suggests the following procedure for proving that the
image in $\R^n$ of certain kinds of maps are automatically convex. Let's
assume that we have some set $E$ and some function $f$ that takes elements
of $E$ into $\R^n$.  To prove that $f(E)$ is convex you can
try the following procedure.  We shall illustrate this in several cases in
Parts II and III of this paper.

\bigskip
\noindent {\bf Automatic Convexity Procedure.}

1.  Find a suitable linear space $X$ and linear map $\Psi$
such that the elements $E$ can be found in $X$ (perhaps in a slightly
different guise) and $f(e) = \Psi(e)$ for each $e \in E$.

2.  Define $K = Conv(E)$ (or perhaps the closure of $Conv(E)$).  Show that the
extreme points of $K$ lie in $E$.

3.  Show that the facial dimension of $K$ is less than the dimension of the
range of $\Psi$, possibly using the intersection theorems in Part I below.

4.  Apply [AA, 1.7] to get the desired convexity.
\bigskip

A knowledge of the facial structure of $K$ is crucial to any
application of [AA, 1.7].  In Part I we prove two new theorems that
describe the facial structure of the intersection of a convex set with
certain subspaces in terms of the facial structure of the original convex
set. These theorems will allow new applications of the automatic convexity
procedure. In Part II of the present paper we discuss
numerical range as an application of pure convexity theory in a way that
(we believe) unravels the mystery and paves the way for more theorems
having convexity as their conclusions.  In Part III we further extend
Lyapunov's convexity theorem and even the continuous lifting theorem of
Samet [S]; again our methods open the way for many more results of the same
type.
\bigskip

\noindent
{\bf PART I:  THE INTERSECTION THEOREMS}
\bigskip

\noindent {\bf Algebraic Intersection Theorem.}  Given a subspace  $H$
in $X$ and a point $x \in X$, let $F$ be a face of $(H+x) \cap K$. Then
$G(K,F) \cap H = F$.
\medskip
\noindent
{\it Proof.}  WLOG we can assume that $x=0$. From [AA, 1.1 and 1.2],
$G(K,v)$ consists of all elements  $y$ of $K$
such that there exists $\lambda > 0$ such that $(1+\lambda)v -
\lambda y\in K$.  Let $G =\bigcup \{G(K,v): v \in F\}$.
Claim $G(K,F) = G$. The inclusion $G \subset G(K,F)$ is clear
from the face property, so we need only show that
$G$ is a face of $K$ and that $G \cap H=F$.

If $x, y \in G$, then there exist $v, w \in F$ such that $x \in
G(K, v)$ and $y \in G(K,w)$.  We can assume a single $\lambda$
such that
$(1+\lambda)v -
\lambda x\in K$ and
$(1+\lambda)w -
\lambda y\in K$.
For any $\alpha \in (0,1)$,
$$
\alpha((1+\lambda)v -\lambda x) + (1-\alpha)((1+\lambda)w -
\lambda y) \in K
$$
Grouping the $(1+\lambda)$ terms and the $\lambda$ terms, we get
$$
(1+\lambda)(\alpha v + (1 - \alpha) w) - \lambda (\alpha x + (1
- \alpha) y)\in K.
$$  Thus $(\alpha x + (1 - \alpha) y) \in G$.
This shows that  $G$ is convex.

To show that  $G$ is a face of $K$, assume $x, y \in K$ such that
$.5(x+y) \in G$.  Then there exists $v \in F$ such that $.5(x+y)
\in G(K,v)$.  But $G(K,v)$ is a face of $K$, so $x, y \in G(K,v)
\subset G $.  Thus $G$ is a face of $K$.

Finally we show that $G \cap H=F$. The inclusion $F \subset G \cap H$
is clear from the definition of $G$. Now if $y \in G \cap H$, then
there is a $v \in F$ such that
$y\in G(K,v)\cap H$. Thus there exists $\lambda > 0$ such that
$(1+\lambda)v-\lambda y \in K$. But $ v \in F \subset H$ and $y \in H$, so
$(1+\lambda)v-\lambda y \in H$ since $H$ is a subspace.  Thus
$(1+\lambda)v-\lambda y \in H \cap K$. Since $F$ is a face of $H \cap K$,
$y\in F$. \hfill\hal

\medskip
\noindent {\bf COMMENT.} If $F$ has a weak internal point $v$ (in the sense
of [AA, p. 8]), then $G(K,F) = G(K,v)$. However,
many interesting infinite dimensional convex sets do not have weak
internal points, e.g.\ the state space of $C([0,1])$ or most any other
interesting C*-algebra.
\medskip

Now we prove a topological version of this result.  As will be clear from
a subsequent example, we need to consider a restricted class of subspaces
$H$ in the topological situation.

\medskip

\noindent {\bf Topological Intersection Theorem:}.  Assume now that $K$
is a convex, closed set in a locally convex space $X$. Given a closed
subspace  $H$ of finite co-dimension in $X$ and a point $y \in X$, let
$F$ be a closed face of $(H+y)\cap K$. Then $G(K,F)$ is closed and
$G(K,F) \cap H = F$.
\medskip

\noindent {\it Proof.} By a simple induction argument, it suffices to prove
the theorem under the asumption that $H$ is a closed hyperplane, and WLOG we
can assume that $y=0$. Let
$f: X \rightarrow \R$ be a continuous linear functional such that $H =
f^{-1}(0)$.  We need only prove that G(K,F) is closed, as
$G(K,F)\cap H = F$ follows from the Algebraic Intersection Theorem.

Suppose $\{x_t\}$ is a net in $G(K,F)$ such that $x_t
\rightarrow x$; we must show $x \in G(K,F)$.  Exchanging $-f$ for $f$ if
necessary and passing to a subnet, we can assume that $f(x_t)
\geq 0$ for all $t$.  If $f(x_t) = 0$ frequently, then we can pass to a
subnet such that each $x_t \in G(K,F) \cap H = F$, and so $x \in F$ (and
hence $x \in G(K,F))$ because $F$ is closed.

Otherwise, pass to a subnet such that $f(x_t) > 0$ for all t.  Let
$x_0$ be any of the $x_t$ and fix it. Since $f(x_0) > 0$, $x_0$ can't lie
in $F$, so by [AA, 1.1]
$x_0 \in G(K,F)$ implies that there is a $y_0 \in G(K,F)$ such that the
open line segment $(x_0, y_0)$ intersects $F$.  It follows by linearity of
$f$ that
$f(y_0) < 0$. Now for each $t$, linearity of $f$ implies that there is a
unique point $z_t$ in $(x_t,y_0)$ such that $f(z_t) = 0$, i.e.\
$z_t \in H\cap (x_t, y_0)$. Explicitly,  $z_t = r_t x_t + (1 - r_t) y_0$
where $r_t = - f(y_0)/(f(x_t) - f(y_0)) \in (0,1)$ since $f(x_0) < 0$.
Since $x_t$ and
$y_0$ are both in $G(K,F)$, it follows from convexity of $G(K,F)$ that $z_t
\in G(K,F)$.  Hence
$z_t\in G(K,F)\cap H = F$ by the Algebraic Intersection Theorem.  Now
$r = \lim r_t = - f(y_0)/(f(x) - f(y_0))> 0$ because
$x_t\rightarrow x$.  Thus
$(z_t)$ converges; let $z = \lim z_t$, so $z
\in F$. Then $z = r x + (1 - r) y_0$, so if $r = 1$ then $x = z \in F$. If
$r < 1$, then the line segment $(x, y_0)$ contains $z$, which implies $x \in
G(K,z) \subset G(K,F)$ . Thus $G(K,F)$ is closed.  \hfill\hal

\bigskip
\noindent {\bf EXAMPLE.} In this example we show why it is necessary to
restrict $H$ to a subspace of finite co-dimension in the Topological
Intersection Theorem.

We work in the Banach space $c_0({\bf Z})$.  Let $h$ be the sequence with
$n$th term
$$h_n = \cases{1/n& if $n \geq 1$\cr 0& if $n \leq 0$\cr}.$$
For each $n \in {\bf Z}$ let $e^n$ be the sequence which is $1$ at $n$ and
$0$ elsewhere.  Let $K_1$ be the closed convex hull of the vectors $k^n =
{1\over n}e^{-n} + h$ (for $n \geq 1$); let $K_2$ be the closed convex
hull of the vectors $-{1\over n}k^n$ (for $n \geq 1$); and let $F$ be the
set of sequences $a = (a_n)$ such that $a_n = 0$ for $n \leq 0$
and $0 \leq a_n \leq n^{-2}$ for $n \geq 1$.

$K_1$ and $K_2$ are each the closed convex hull of a convergent sequence
of vectors in a Banach space, and hence are compact. $F$ is compact
because it is closed and totally bounded. Thus the convex hull $K$ of
$K_1$, $K_2$, and $F$ is compact. (It is a continuous image of the
compact set $K_1 \times K_2 \times F \times S$ where $S = \{(r,s,t):
r, s, t \geq 0$ and $r + s + t = 1\}$.)

Explicitly, $K_1$ is the set of sequences $(a_n)$ such that $a_0 = 0, a_n
= h_n$ for $n \geq 1$, $a_n \geq 0$ for $n \leq 0$, and
$\sum_{n = 0}^\infty na_{-n} \leq 1$. $K_2$ is the set of sequences
$b_n$ such that $b_0 = 0, b_n \leq 0$ for $n \leq 0$, $\sum_{n=0}^\infty
n^2b_{-n} \geq -1$, and $b_n = \alpha h_n$ for $n \geq 1$ where
$\alpha = \sum_{n=0}^\infty nb_{-n}$.

Observe that if $ra + sb$ is in the convex hull of $K_1$ and $K_2$
($a \in K_1$, $b \in K_2$, $r + s = 1$) and $ra_n + sb_n = 0$ for all
$n \leq 0$, then we must have
$$s\cdot \sum_{n = 0}^\infty nb_{-n} = -r\cdot \sum_{n = 0}^\infty na_{-n}
\geq -r.$$
Thus $ra + sb = rh + s\alpha h$ where $\alpha \geq -r/s$, and thus
$ra + sb = \beta h$ with $0 \leq \beta \leq 1$.

Now let $H$ be the set of sequences $(a_n)$ such that $a_n = 0$ for all
$n \leq 0$. This is a closed subspace of $c_0({\bf Z})$. $K_1$
intersects $H$ in the point $h$ and $F$ is contained in $H$, so
$K \cap H$ contains the convex hull $C$ of $h$ and $F$. Moreover,
any element of $K$ --- that is, any convex combination $ra + sb + tc$
with $a \in K_1$, $b \in K_2$, and $c \in F$ --- which lies in $H$
must satisfy $ra + sb \in H$; then by the last paragraph,
$ra + sb = \beta h$ where $0 \leq \beta \leq r + s$, so we have
$$ra + sb + tc = \beta h + tc = (1-t){{\beta}\over{1-t}}h + tc$$
where $\beta/(1-t) = \beta/(r + s) \leq 1$. Since $C$ contains $h$ and $0$,
it contains $(\beta/(1-t)) h$, and therefore it contains $ra + sb + tc$.
We have shown that $C = K \cap H$.

Next we claim $F$ is a closed face of $C$. It is closed because it
is compact.  It is a face because if $a, b \in C$ and neither belongs
to $F$ then $\lim a_n/h_n$ and $\lim b_n/h_n$ both exist and are
strictly positive, so the same is true of $(a + b)/2$, which implies
$(a + b)/2 \not\in F$. This proves the claim.

Finally, we claim that any closed face $G$ of $K$ that contains $F$
must contain $h$. For $0 \in F$, and $0$ lies in the line segment
joining $k^n$ and $-{1\over n}k^n$, which both belong to $K$, so
$k^n \in G$. Since $G$ is closed and $k^n \rightarrow h$, it
follows that $h \in G$. This proves the final claim and shows that
$G \cap H \neq F$.  \hfill\hal
\bigskip

\medskip
\noindent {\bf PART II:  APPLICATIONS TO NUMERICAL RANGE}
\medskip

\noindent {\bf NOTATION:}  Let $M_n$ denote the set of $n \times n$ complex
matrices and $U_n$ the set of unitary matrices in $M_n$.  Let
$\tau$ denote the trace on
$M_n$ and $1$ the identity matrix. For $a, b \in M_n$ write $a \ge b$ if
$a-b$ is positive semi-definite. Define $K=\{a\in M_n : 0 \le a \le 1\}$. When
we need to specify a norm on $M_n$ we shall always take the operator 
norm, i.e.\
$\|a\| = sup\{\|a\eta\|: \eta$ is a unit vector in $\R^n\}.$  The $k$-numerical
range of an $n \times n$ matrix $b$ is
$W_k(b) =\{(1/k)\sum_1^k (bx_j,x_j):$ the $x_j$ are orthonormal$\}$.
When no confusion can develop we identify the complex plane with $\R^2$.
\medskip

Let's illustrate our four step method by proving the convexity of the
$k$-numerical range of $b\in M_n$.  This was first shown by Berger [B]. A more
accessible proof based on the convexity of the ordinary numerical range can be
found in [Hal-2, Problem 167].  The first step is to linearize the function
that produces the points in the $k$-numerical range.  The definition 
of $W_k(b)$
calls for calculating a complex number for each $k$-tuple of 
orthonormal vectors
in $\R^n$.  Replace such a $k$-tuple $\{x_j\}_1^k$ with the 
orthogonal projection
$p$ of their span.  Then $\tau(pb) = \sum_1^k (bx_j,x_j)$, thus we can see
that $kW_k(b) =\{\tau(bq) : q$ is a projection of rank $k\}$, so it suffices
to show that the latter set is convex.  Setting $E=\{q \in M_n: q$ is a
projection of rank $k\}$ completes the first step.

For the second step we define $Q_k= Conv(E)$.  Clearly $Q_k \subset K$.  Since
$\tau(q) = k$ for each $ q\in E$, this suggests that we consider $Q_k$ as a
subset of $\{a \in K : \tau(a) = k\}$.  In Proposition 1 below we show
that $E(Q_k) = E$.

For the third step we need to determine the facial structure of $Q_k$.  It
seems
sensible to start with $K$.  This is probably classical, but a readable (and
more general) account appears as [AP, 2.2] where faces of $K$ are shown to have
the form $F=\{x\in K : p \le x \le q\}$,
where $q,p$ are self-adjoint projections in $M_n$. This can be rewritten
in terms of the difference $q-p = r$ as $F=p+rKr$. A face of this form is an
extreme point exactly when $r=0$, and then the extreme point is just the
projection $p$. i.e.\ the extreme points of $K$ are exactly the projections.
Since the analysis of the facial structure of $Q_k$ uses the intersection
theorem from Part I, we state the facts as a proposition.

\medskip

\noindent {\bf PROPOSITION 1.} For $1 \le k < n$,
$Q_k=\{a \in K : \tau (a) = k\}$. The facial dimension of $Q_k$ is 3.
Further, the extreme points of $Q_k$ are exactly the projections of rank $k$.
\medskip

\noindent Proof:  We already noted that
$Q_k \subset \{a \in K : \tau (a) = k\}$. If we show that the right hand side
has exactly the projections of rank $k$ as its extreme point set, then equality
will follow.

Note that if we intersect $K$ with the hyperplane $H=\{x : \tau(x) = k\}$, then
we get exactly $\{a \in K : \tau (a) = k\}$. Using the notation developed just
above the statement of the proposition, let a face $F$ of $K$ have the form
$F=p+rKr$. By the Algebraic Intersection Theorem the typical face of
$\{a \in K : \tau (a) = k\}$ is $F \cap H$.

If $\tau(r) = 2$, then the face $F$ has real dimension  4 since
this is easily verified for $rKr$. If $\tau(r) > 2$, then
the dimension of $F\cap H$ is even larger. On the
other hand, if $rank(q-p) = 1$, then $F$ is exactly the line segment joining
$p$ and $q$.  Such a line segment can meet $H$ only at one of the end
points, i.e.\ at a projection of rank $k$. Thus we have shown that the set of
extreme points of
$\{a\in K : \tau (a) = k\}$ is exactly the set of projections of rank k,
thereby completing the proof of $Q_k=\{a \in K : \tau (a) = k\}$.  We also
have shown that $Q_k$ has no faces of dimension 1 or 2, hence its facial
dimension is at least 3. Faces of dimension exactly 3 occur when $\tau(r)
= 2$.\hfill\hal
\bigskip

We complete step 4 with the following proposition.
\medskip

\noindent {\bf PROPOSITION 2.} If $b \in M_n$, then the $k$-numerical range
of $b$ is convex.
\medskip
\noindent {\it Proof.} The linear map $\Psi(a) = \tau(ab)$ takes $Q_k$ into
$\C$ and its range is exactly the $k$-numerical range of $b$.  Since the
facial dimension of $Q_k$ is 3 and the extreme points are projections of
rank $k$, [AA, 1.7] gives the desired convexity.  \hfill\hal
\bigskip

As another example of this method, we prove the convexity of the
$c$-numerical range for a self-adjoint element c of $M_n$.
   For any $c \in M_n$ the $c$-numerical
range of a matrix $a \in M_n$ is defined to be
$W_c(a) =\{\tau (cu^*au), u \in U_n\}$. It is easy to check that the
$k$-numerical range $W(a)$ is obtained from this definition when
$c$ is taken to be a self-adjoint projection matrix of rank $k$ (for the
ordinary numerical range simply take $k = 1$).   It is known that the
$c$-numerical range is convex when
$c$ is self-adjoint [GR, sect. 5.5]. In the next proposition we show
``why'' this is true.
\medskip

\noindent {\bf PROPOSITION 3.} If $c \in M_n$ is self adjoint, then the
$c$-numerical range of $b$ is convex for all $b \in M_n$.

\medskip
\noindent {\it Proof.} In this formulation the first step of the automatic
convexity procedure is straightforward.  Fix an element $b \in M_n$. Define
$E=\{u^*cu : u\in U_n\}$.  Let $M = Conv(E)$. Note that $M$ is closed since
$E$ is closed [W, 2.2.6]. The set of extreme  points of $M$
is exactly the set $E$ since
$M$ contains extreme points [W, 2.6.16], these lie in $E$ [W, 2.6.4], and any
point of $E$ can be mapped onto any other by a linear isometry of $M$ onto
itself (namely $u^*cu \rightarrow v^*u(u^*cu)u^*v = v^*cv$). For any
$a\in M_n$ define $\Psi(a) = \tau(ab)$.  Then $\Psi(E)$ is exactly the
$c$-numerical range of $b$.  To complete the proof using [AA, 1.7], 
we need only
show that the facial dimension of $M$ is at least 3.  This is done in the
following lemma.  \hfill\hal

\bigskip

  Since we have to borrow from matrix theory
for the proof of the next lemma, for comparison and convenience we
use the notation of [GR, Section 5.5]. Because of the change to the
notation of [GR, Section 5.5] what we called $c$ in the previous
proposition is now $C$, while $c$ stands for the real vector
consisting of the eigenvalues of $C$.
\bigskip
\noindent {\bf LEMMA 3.5.} Fix a self-adjoint matrix $C \in M_n$. Then the
facial dimension of
$M=\{U^*CU : U \in M_n\}$ is at least 3.
\medskip
   \noindent {\it Proof.}
Let $\alpha, \beta \in \R^k$. We say that $\alpha$ is obtained from
$\beta$ by pinching if all components of $\alpha$ and $\beta$ agree
except for two, $\alpha_i$ and $\alpha_j$, which satisfy
$\alpha_i = \lambda\beta_i + (1-\lambda)\beta_j$ and
$\alpha_j = (1-\lambda)\beta_i + \lambda\beta_j$ for some $\lambda
\in [0,1]$. We require the following fact: the positive vector
$\alpha$ is obtained from the positive vector $\beta$ by a finite number of
pinchings if and only if $$\sum_{i=1}^k \alpha_i \leq \sum_{i=1}^k \beta_i$$
for $1 \leq k \leq n$, with equality when $k = n$. Write $\alpha \prec
\beta$ for this relation.

   Since adding a scalar multiple of
the identity matrix to $C$ only shifts the $C$-numerical range, WLOG we can let
$C$ be the positive diagonal matrix with diagonal $c$, denoted $C = [c]$, where
$c$ is arranged in decreasing order. Let $M' = \{U^*[b]U: U$ is unitary and $b
\prec c\}$. We shall show that $M' = M$. Note that $M'$ is the set of positive
matrices
$B$ whose ordered eigenvalue list $b$ satisfies $b \prec c$. Observe that the
sum of the first
$k$ eigenvalues of $B$ equals $\sup \{ \tau(BP): P$ is a rank $k$
projection$\}$ [AAW, Lemma 1.3]. Thus, $M'$ is the set of positive
matrices $B$ such that $\tau(B) = \tau(C)$ and
$$\tau(BP) \leq \sum_{i=1}^k c_i$$
for $1 \leq k \leq n$ and every rank $k$ projection $P$. It easily follows
that $M'$ is closed and convex.

Next, we claim that the extreme points of $M'$ are precisely the matrices
of the form $U^*[c]U$ for $U$ a unitary matrix. To see this, let $B =
U^*[b]U \in M'$ and suppose $B$ is not of the form $U^*[c]U$. Then $[b]$
is obtained from $[c]$ by a finite, nonempty sequence of pinchings. It
follows that $[b]$ is obtained from some $[a]$ by a single pinching, where
$a \prec c$. That is, $b_i = ta_i + (1-t)a_j$ and $b_j = (1-t)a_i + ta_j$
for some $t \in (0,1)$, where $a_i \neq a_j$, and all other components of
$a$ and $b$ agree. Let $a'$ be the real vector obtained from $a$ by switching
the $i$ and $j$ components. Then $A = U^*[a]U$ and $A' = U^*[a']U$ are both
in $M'$, and $B = tA + (1-t)A'$. So $B$ is not an extreme point. Thus, every
extreme point of $M'$ must be of the form $U^*[c]U$.  Thus $M'=M$ by [W,
2.6.16].

Finally, we claim that the facial dimension of $M$ is at least 3. To see
this, let $B = U^*[b]U \in M$ and suppose $B$ is not an extreme point.
Define $A$ and $A'$ as in the last paragraph. Then
$$A_{[ij]} = \left[\matrix{a_i&0\cr 0&a_j\cr}\right]\qquad{\rm and}
\qquad A'_{[ij]} = \left[\matrix{a_j&0\cr 0&a_i\cr}\right],$$
where we use the subscript $[ij]$ to indicate restriction to the
$(i,i)$, $(i,j)$, $(j,i)$, and $(j,j)$ entries. (Recall that $A$ and
$A'$ agree elsewhere.) Define new matrices $A_1$, $A_1'$, $A_2$, and
$A_2'$ by setting
$$\eqalign{(A_1)_{[ij]}
&= \left[\matrix{ta_i + (1-t)a_j&(t-t^2)^{1/2}(a_i-a_j)\cr
(t-t^2)^{1/2}(a_i-a_j)&(1-t)a_i + ta_j\cr}\right]\cr
(A_1')_{[ij]}
&= \left[\matrix{ta_i + (1-t)a_j&-(t-t^2)^{1/2}(a_i-a_j)\cr
-(t-t^2)^{1/2}(a_i-a_j)&(1-t)a_i + ta_j\cr}\right]\cr
(A_2)_{[ij]}
&= \left[\matrix{ta_i + (1-t)a_j&i(t-t^2)^{1/2}(a_i-a_j)\cr
-i(t-t^2)^{1/2}(a_i-a_j)&(1-t)a_i + ta_j\cr}\right]\cr
(A_2')_{[ij]}
&=  \left[\matrix{ta_i + (1-t)a_j&-i(t-t^2)^{1/2}(a_i-a_j)\cr
i(t-t^2)^{1/2}(a_i-a_j)&(1-t)a_i + ta_j\cr}\right]\cr}$$
and letting them agree with $A$ and $A'$ elsewhere. It is clear that
each of these matrices is self-adjoint, and as the $2\times 2$ parts
all have the same trace and determinant, they all have the same eigenvalues
(namely, $a_i$ and $a_j$). Thus they all belong to $M$. But $B =
(A_1 + A_1')/2 = (A_2 + A_2')/2$, and the affine space spanned by $A$,
$A'$, $A_1$, $A_1'$, $A_2$, and $A_2'$ is three-dimensional, so the smallest
face containing $B$ has dimension at least 3. This proves the final claim.
   \hfill\hal

\bigskip
\noindent
{\bf PART III:  APPLICATIONS TO LYAPUNOV TYPE THEOREMS}
\medskip
Let $(X,\M)$ be a measurable space. A vector measure is an $n$-tuple
$(\mu _1, \ldots, \mu_n) = {\bf \mu}$ of real-valued measures on $(X, \M)$.
Lyapunov's Theorem [L] states that the range of $\mu$ is a convex, compact
set in $\R^n$.  Following the 4 step plan for proving convexity (and often
compactness in the same stroke, as is the case here) one observes that
$$\int \chi_A d\mu = \mu(A) = (\int \chi_A d\mu_1, \ldots, \int \chi_A
d\mu_n) = (\int \chi_A f_1 d\nu, \ldots, \int \chi_A
   f_n d\nu),$$ where $\chi_A$ is the characteristic function of the set
$A$, $\nu = \sum^n_1 |\mu_i|$ is a finite, positive measure, and $f_i$ is
the Radon-Nikodym derivative of $\mu_i$ with respect to $\nu$ for each
$i$.  This formulation suggested the definition of the map $\Psi :
L^{\infty}(X, \M, \nu)
\rightarrow \R^n$ by $\Psi(g) = (\int g f_1 d\nu, \ldots, \int g f_n d\nu).$
Moving to step 2 in the plan, we note that if $E$ is viewed as the set of
characteristic functions in $L^\infty(X, \M, \nu)$, then the closed convex
hull $K$ of
$E$ in the weak* topology is exactly the set of positive functions of norm
no more than 1. The facial dimension of $K$ is shown in [AP] to be
$\infty$, so [AA, 1.7], the weak* compactness of $K$ and the weak*
continuity of $\Psi$ complete the proof of Lyapunov's Theorem.

As with the numerical range situation discussed earlier in this paper, once
the problem was put into the correct notation, the convexity was automatic
from facial structure considerations and [AA, 1.7].  Of course [AA]
contained many results that could be viewed as generalizations of
Lyapunov's Theorem.  Now let's combine these results with the Topological
Intersection Theorem to show how even more theorems of the Lyapunov type
are true using our 4 step method. In the next Theorem we extend [AA, 2.5],
which is itself an extension of Lyapunov's theorem to a non-commutative
situation.

\medskip
\noindent {\bf THEOREM 4.} Suppose that $N$ is a non-atomic von Neuman
algebra and $\{f_1, \ldots, f_n\}$ and $\{g_1, \ldots, g_k\}$ are self-adjoint,
normal linear functionals on $N$.  Let $z_1, \ldots, z_n \in \R$ and define
$$K=\{a \in N : \|a\| \le 1, a \ge 0, f_j(a) = z_j, j= 1, \ldots, n\}.$$
Let $N_{sa}$ denote the set of self adjoint elements of $N$.  Define
$\Psi : N_{sa} \rightarrow\R^k$ by
$\Psi(a) = (g_1(a), \ldots, g_k(a))$.  Then $E(K)= \{p :$ $p$ is a
projection in $K\}$ and $\Psi(K)=\Psi(E(K)).$

If $N$ is abelian, then there is a continuous map $\Phi : \Psi(K) \rightarrow
E(K)$ that is a right inverse for $\Psi$.
\medskip

\noindent {\it Proof.} If $K$ is void, the theorem is trivially true, so
asume not. If
$N^+_1$ denotes the positive part of the unit ball of $N$, then the facial
dimension of
$N^+_1$ is
$\infty$ by [AP, 2.2]. Since
$K$ is the intersection of
$N^+_1$ with a subspace of finite codimension, the Topological
Intersection Theorem applies to show that the faces of $K$ are either
extreme points of
$N^+_1$ or else infinite dimensional faces.  Since the extreme points of
$N^+_1$ are exactly the projections of $N$ by [AP, 2.2] , we get 
$E(K) =\{p :$ $p$ is a
projection in $K\}$. The conclusion $\Psi(K)=\Psi(E(K))$
follows from [AA, 1.7].

Now assume that $N$ is abelian.  Define $\Psi' : N_{sa} \rightarrow \R^{n+k}$
by the formula
$$\Psi'(a) = (\Psi(a), f_1(a), \ldots, f_n(a)).$$
By Lyapunov's Theorem $\Psi'(N^+_1)$ is compact and convex.  By [S] there is a
continuous right inverse $\Phi' : \Psi(K) \rightarrow
E(K)$ for $\Psi'$.  Now let $$S=\{(\Psi(a), f_1(a), \ldots, f_n(a)) \in
\Psi'(N^+_1) : f_j(a) = z_j, j = 1, 2, \ldots, n\}.$$
Clearly $\Psi'^{-1}(S) = \{a \in N^+_1 : f_j(a) = z_j, j = 1, 2, \ldots, n\}
= K$. Thus the restriction of $\Phi'$ to $K$ is the desired lifting if we
identify the first $k$ coordinates of $\R^{k+n}$ with $\R^k$.
\hfill\hal
\bigskip

We present two corollaries of Theorem 4. The first is
a version of Lyapunov's theorem with linear constraints, and could possibly
have applications in control theory along the lines of the classical Lyapunov
theorem [HLS]. The second gives a von Neumann algebra version of the
convexity of the $k$-numerical range (where here $k = z$).
\bigskip

\noindent {\bf COROLLARY 5.} Let $(X, \M)$ be a measurable space, let
$\mu = (\mu_1, \ldots, \mu_k)$ be a vector measure on $(X,\M)$, let
$\nu_1, \ldots, \nu_n$ be measures which are absolutely continuous with
respect to $\nu = |\mu_1| + \cdots + |\mu_k|$, and let $z_1, \ldots, z_n \in
\R$. Then the set $\{\mu(A): \nu_j(A) = z_j, j = 1, 2, \ldots, n\}$ is
compact and convex.
\medskip

\noindent {\it Proof.} We translate into the language of Theorem 4 by
letting $N = L^\infty(X, \nu)$ and letting the $f_i$ and $g_j$ be the
Radon-Nikodym derivatives of the $\mu_i$ and $\nu_j$ with respect to
$\nu$. Then $\{\mu(A): \nu_j(A) = z_j, j = 1, 2, \ldots, n\} =
\Psi(E(K)) = \Psi(K)$, and $\Psi(K)$ is clearly compact and convex.\hfill\hal
\bigskip

\noindent {\bf COROLLARY 6.} Let $N$ be a non-atomic von Neumann algebra
with normal tracial state $\tau$ and let $b \in N$ and $z \in [0,1]$. Then
the set
$$\{\tau(pb): p\hbox{ is a projection and }\tau(p) = z\}$$
is a compact and convex subset of $\C$.
\medskip

\noindent {\it Proof.} In Theorem 4, take $n = k = 1$, $f_1 = \tau$,
$g_1 = \tau(\cdot b)$, and $z_1 = z$.\hfill\hal

\bigskip
\centerline {\bf References}
\medskip
\noindent [AA] C.\ Akemann and J.\ Anderson, {\it Lyapunov
Theorems for Operator Algebras}, Memoirs of the Amer. Math. Soc. No. 458,
v. 94, Providence, November 1991.
\medskip

\noindent [AAW] C.\ Akemann, J.\ Anderson, and N.\ Weaver, A geometric
spectral theory for $n$-tuples of self-adjoint operators in finite von
Neumann algebras, {\it J. Funct. Anal.} 165 (1999), 258-292.
\medskip

\noindent [AP] C.\ Akemann and G.\ K.\  Pedersen,  Facial
structure in operator algebra theory, {\it Proc. London Math. Soc.} (3) 64
(1992) 418-448.
\medskip

\noindent [ASW] C.\ Akemann, G.\ Shell, and Nik \ Weaver, Locally
nonconical convexity, {it J. Convex Analysis}, 8 (2001), 1-21.
\medskip

\noindent [B] C.\ Berger, Normal Dilations, PhD. Dissertation, Cornell, 1963.
\medskip

\noindent [GR] K.\ Gustafson and D.\ Rao, {\it Numerical Range: The
filed of values of operators and matrices}, Springer, New York, 1997.
\medskip

\noindent [Hal-1] P.\ Halmos, On the set of values of a finite measure,
{\it Bull. Amer. Math. Soc.} 53, (1947)
\medskip

\noindent [Hal-2] P.\ Halmos, {\it A Hilbert Space Problem Book}, Van
Nostrand, New York, 1967.
\medskip

\noindent [Hau] F.\ Hausdorff, Das algebraische Analogon zu einem Satz von
Fej\'er, {\it Math. Z} {\bf 2} (1918), 187-197.
\medskip

\noindent [HLS] H.\ Hermes and J.\ P.\ LaSalle, {\it Functional analysis
and time optimal control}, Academic Press, New York, 1969.
\medskip

\noindent [HJ] R.\ Horn and C.\ Johnson, {\it Topics in Matrix
Analysis}, Cambridge University Press, Cambridge, 1991.
\medskip

\noindent [Li] J.\ Lindennstrauss, A short proof of Lyapunov's convexity
theorem, {\it J. Math. Mech.} 15 (1966), 971-972.
\medskip

\noindent [Ly] A.\  Lyapunov, On completely additive vector functions,
Bull. Akad. Sci. USSR 4 (1940), 465-478. (Russian)
\medskip

\noindent [P] Y.\ Poon, Generalized numerical ranges, joint positive
definiteness and multiple eigenvalues, {\it Proc. Amer. Math. Soc.}
(6) {\bf 125} (1997), 1625-1634.
\medskip

\noindent [S]  D.\ Samet, Continuous selections for vector measures,
{\it Mathematics of operations Research}, {\bf 12}, (3) (1987), 536-543.
\medskip

\noindent [T] O.\ Toeplitz, Der Wertvorrat einer Bilinearform, {\it Math.
Z} {\bf 3} (1919), 314-316.
\medskip

\noindent [W] R.\ Webster, {\it Convexity}, Oxford University Press, Oxford,
1994.

\bigskip
\bigskip
Charles A. Akemann

Department of Mathematics

University of California

Santa Barbara, CA 93106, USA

akemann@math.ucsb.edu
\bigskip
Nik Weaver

Department of Mathematics

Washington University

St. Louis, MO 63130, USA

nweaver@math.wustl.edu

\end